\documentclass[a4paper,11pt]{article}

\usepackage{latexsym}
\usepackage{amsmath,amsfonts,amssymb,amsthm}
\newtheorem{thm}{Theorem}[section]
\newtheorem{lem}[thm]{Lemma}

\newtheorem{cor}[thm]{Corollary}
\newtheorem{prop}[thm]{Proposition}
\newtheorem{conj}[thm]{Conjecture}
\newtheorem{rem}[thm]{Remark}

\input epsf

\title{Fair Triangulations}
\author{Roland Bacher}
\begin{document}
\maketitle
%\par checkerbrd1.tex dans recherche/catalanwinding

{\sl Abstract: We describe the statistics of checkerboard triangulations
obtained by colouring black every other triangle in triangulations of
convex polygons.}\footnote{Math. Class: 05A15, Keywords: Catalan number, checkerboard triangulation, dissections of polygons}

\section{Introduction}

Triangulations of strictly convex polygons are classical combinatorial
objects studied almost ever since the sad moment when 
mankind made his biggest mistake
and left his favourite tree. Out of remorse and longing for the lost
paradise, mankind (in fact, only a small part of it, those successful 
with girls had more urgent and important matters to attend to)
started the study of trees and discovered that
plane rooted trees and triangulated convex polygons share many combinatorial
properties. This caused some fellows to investigate triangles 
thus giving birth to geometry and eventually analysis and real numbers.
Other chaps, tempted by enumeration, discovered the natural numbers 
which led later to arithmetics and algebra. 
They tried successfully to enumerate
plane rooted trees (or, equivalently, triangulations of convex polygons) 
by introducing the now famous sequence 
of the ubiquituous Catalan numbers, namely
$$1,1,2,5,14,42,132,429,1430,4862,16796,58786,208012,742900,\dots,$$ 
see sequence A108 in \cite{EIS}. Exercice 6.19 in 
\cite{St2} lists $66$ sequences of sets enumerated by them.
\cite{St2} contains also some historical information
which is obviously rather less reliable than the translation of
the Holy Scriptures presented above.

The generating series of the Catalan sequence is the algebraic function
$$g(x)=\frac{1-\sqrt{1-4x}}{2x}=\sum_{n=0}^\infty {2n\choose n}
\frac{x^n}{n+1}\ .$$
It satisfies the equation $g(x)=1+xg(x)^2$.

This paper is devoted to the study of triangulations which are
coloured like a chessboard: every triangle is either black or white 
and edge-adjacent triangles have different colours. Such a colouring
is unique, up to permutation of the colours.
Call a triangulation $\tau$ {\it fair} if its two colourings
involve the same number of black and white triangles.

\begin{thm} \label{thmfair} Every strictly convex polygon with $2n+2$
vertices has exactly ${2n\choose n}^2\frac{1}{n+1}$
fair triangulations.
\end{thm}

The choice of a marked edge $e_*$ on the boundary $\partial P$
of a polygon $P$ having at least $3$ vertices
selects in any triangulation of $P$ a marked triangle $\Delta_*$, 
which contains the marked edge $e_*$. The imposition of 
the colour black for the marked triangle $\Delta_*$
removes the colour-ambiguity. We get in this way a bijection between
triangulations and ``checkerboard triangulations'' of 
convex polygons $P\supset e_*$ having a marked edge.
Section \ref{sectchecktriang} contains our main result,
Theorem \ref{thmwindcoeffs}, which describes closed formulae
for the the number of such checkerboard triangulations
involving $n_b$ black and $n_w$ white triangles. 
Theorem \ref{thmfair} corresponds to the special case $n_b=n_w$.
The proof of Theorem \ref{thmwindcoeffs} uses probably well-known
techniques usable for similar problems: 
algebraic generating series which can be expressed as finite
sums of a kind of hypergeometric series in several variables.

We study also $d-$dissections, a generalisation
of triangulations corresponding to the special case $d=2$. 
(The choice of the notation
can be motivated by higher dimensional analogues.) A 
$d-$dissection of a strictly convex polygon $P$ is a decomposition 
of $P$ into polygons with $(d+1)-$vertices,
all contained in the vertex-set of $P$.  
Defining fair $d-$dissections in the obvious way,
we have the following conjecture:

\begin{conj} \label{conjdfair} A strictly convex polygon with $2(d-1)n+2$
vertices has exactly ${dn\choose n}^2\frac{1}{(d-1)n+1}$
fair $d-$dissections.
\end{conj}

Conjecture \ref{conjdfair} holds of course for $d=2$ by Theorem 
\ref{thmfair}. For general $d$ it is again a special case: 
Conjecture \ref{conjcoeffwd} of Section \ref{sectmainddiss}
gives  the number of all $d-$dissections of a convex polygon $P$
involving exactly $n_b$ black and $n_w$ white polygons such that 
a prescribed marked edge $e_*\subset P$ belongs to a black triangle.
Obvious modifications of the proof
of Theorem \ref{thmwindcoeffs} establish Conjecture \ref{conjcoeffwd} 
for a few small values of $d$ and show the truth of
Conjecture \ref{conjdfair} for $d\leq 6$.

\begin{rem} The sequence $1,2,12,100,980,10584,\dots,
{2n\choose n}^2\frac{1}{n+1},\dots$ enumerating fair 
triangulations appears as A888 in \cite{EIS},
together with the following description, due to D. Callan: ${2n\choose n}^2
\frac{1}{n+1}$ is the number of lattice paths consisting of $2n$ 
steps in $\{(\pm 1,0),(0,\pm 1)\}$
starting at $(0,0)$ and ending on the diagonal $x=y$
with the constraint of remaining in the tilted quarter plane
$-x\leq y\leq x$. A short bijective proof by Callan (\cite{Ca})
is as follows: 
The pair of orthogonal projections onto the lines
$y=-x$ and $y=x$ induces a bijection between such paths and the 
product of Dyck paths of length $2n$ and of positive paths of 
length $2n$. Dyck paths of length $2n$, given by lattice walks of $2n$ 
steps $\pm 1$ on $\mathbb N$, starting and ending 
at the origin, are enumerated by the Catalan number ${2n\choose n}
\frac{1}{n+1}$. Positive paths of length $2n$,
given by lattice walks of $2n$ steps $\pm 1$ on $\mathbb N$,
starting at the origin, are enumerated by the central binomial
coefficient ${2n\choose n}$.
\end{rem}

\begin{rem} A separate paper will deal with vertex-colourings of
triangulations. The vertices of a $\tau-$triangulated polygon 
$P$ can be coloured with $3$ colours such that every triangle of 
$\tau$ has vertices of all $3$ colours. Such a $3-$colouring
is unique up to permutations of the three colours. 
One can show that a strictly convex polygon
with $3n$ vertices has exactly ${2n-2\choose n-1}^3\frac{3n-2}{n^2}$
triangulations which involve each colour a common number $n$ of times
in every such $3-$colouring of its vertices.
Methods and techniques are completely analogous to those
used in the present paper (although the formulae are simpler
and the necessary symbolic computations heavier). 
\end{rem}

The next two sections contain definitions and our main results.
The rest of the paper is devoted to proofs and a few complements. 

\section{Main results for checkerboard triangulations}\label{sectchecktriang}

Through the rest of this paper, $P_n$ denotes
always some fixed strictly convex polygon with $n+2$ vertices and edges.
We write $P_n\supset e_*$ if $P_n$ is decorated with a marked 
edge $e_*\subset \partial P_n$. 
A {\it triangulation} of $P_n$ is a decomposition of $P_n$ into $n$
non-overlapping triangles. A triangulation will always
be denoted by the greek letter $\tau$.
A triangulation $\tau$ 
of $P_n\supset e_*$ selects a unique marked triangle $\Delta_*\in\tau$ 
of $\tau$ such that $\Delta_*$ contains the marked edge $e_*$. 
We will generally omit
a separate discussion of the degenerate and trivial initial case
$n=0$ corresponding to a polygon reduced to a (double) edge having 
a unique triangulation involving no triangles.

Given a triangulation $\tau$ of $P_n\supset e_*$,
the associated {\it checkerboard colouring} 
partitions the set of all $n=n_b(\tau)+n_w(\tau)$ 
triangles of $\tau$ into a subset 
of $n_b=n_b(\tau)$ black triangles containing the marked
triangle $\Delta_*$ and a complementary subset of
$n_w=n_w(\tau)$ white triangles 
such that edge-adjacent triangles have different colours. 
A {\it checkerboard triangulation} 
is a triangulation $\tau$ of $P_n\supset e_*$
whose $n$ triangles are coloured by the unique checkerboard colouring
associated to $\tau$.

\begin{figure}[h]\label{fig1}
\epsfysize=1.4cm
\centerline{\epsfbox{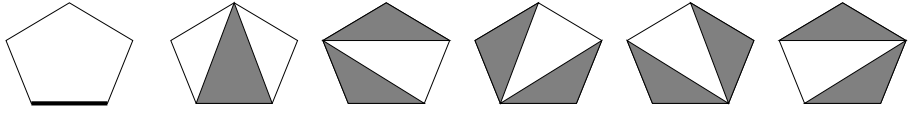}}
\caption{$P_3\supset e_*$ with $e_*$ given by the fat edge and
its $C_3=5$ checkerboard triangulations.}
\end{figure}

Endowing a checkerboard triangulation $\tau$ of $P_n\supset e_*$ 
with the weight $t^{n_b(\tau)-n_w(\tau)}$ and summing the monomials
$t^{n_b(\tau)-n_w(\tau)}$ over
the set $\mathcal T_n$ of all ${2n\choose n}\frac{1}{n+1}$ 
checkerboard triangulations of $P_n\supset e_*$ we get Laurent-polynomials 
$$\begin{array}{l}
\displaystyle w_0(t)=1\\
\displaystyle w_1(t)=t\\
\displaystyle w_2(t)=2\\
\displaystyle w_3(t)=t^{-1}+4t\\
\displaystyle w_4(t)=12+2t^2\\
\displaystyle w_5(t)=12t^{-1}+30t\\
\displaystyle w_6(t)=4t^{-2}+100+28t^2\\
\displaystyle w_7(t)=140t^{-1}+280t+9t^3\\
\displaystyle w_8(t)=90t^{-2}+980+360t^2\\
\displaystyle w_9(t)=22t^{-3}+1680t^{-1}+2940t+220t^3\\
\displaystyle w_{10}(t)=1540t^{-2}+10584+4620t^2+52t^4\\
\quad\vdots\\
w_n(t)=\sum_{\tau\in\mathcal T_n} t^{n_b(\tau)-n_w(\tau)}\in\mathbb N
[t,t^{-1}]
\end{array}$$
recursively defined by the following easy result.

\begin{thm} \label{thmeqW} (i) We have $w_0(t)=1$ and 
$$w_{n+1}(t)=t\sum_{k=0}^nw_k(t^{-1})w_{n-k}(t^{-1})\ .$$

\ \ (ii) The generating series 
$W=\sum_{n=0}^\infty w_n(t)x^n\in\mathbb N[[t,t^{-1},x]]$
satisfies the algebraic equation
$$t(1+tx)-tW+2tx^2W^2+x^3W^4=0\ .$$ 
\end{thm}

Observe that we have $W\in\mathbf R$ with $\mathbf R$ 
denoting through the rest of the paper the 
algebra $\mathbf R=\left(\mathbb Q [t,t^{-1}]\right)[[x]]$
of formal power series in $x$ with coefficients given by 
Laurent-polynomials in $t$.

We have moreover the following result for the Laurent-polynomials
$w_n=w_n(t)$ and the associated generating series $W=\sum_{n=0}^\infty
w_nx^n\in\mathbf R$.

\begin{thm} \label{thmcoeffid} We have the identity
$$(n+2+3j)(w_n,t^{-j})=(n+2-3j)(w_n,t^j)$$
where $(w_n,t^j)$ denotes the coefficient of $t^j$ in 
$w_n=\sum_{j\in\mathbb Z} (w_n,t^j)t^j\in\mathbb N[t,t^{-1}]$.
\end{thm}

Setting ${\overline W}=\sum_{n=0}^\infty w_n(t^{-1})x^n\in\mathbf R$, 
Theorem \ref{thmcoeffid}
amounts to the equality
$$x{\overline W}_x+2{\overline W}+3t{\overline W}_t=xW_x+2W-3tW_t$$
where $W_x=\frac{\partial }{\partial x}W,\ 
W_t=\frac{\partial }{\partial t}W,\ 
{\overline W}_x=\frac{\partial }{\partial x}{\overline W},\ 
{\overline W}_t=\frac{\partial }{\partial t}{\overline W}$.

\begin{cor} \label{cordiffeq}
The generating series $W=\sum_{n=0}^\infty w_nx^n$
satisfies the partial differential equation
$$2t(1-W)+x(2xW-t)W_x+3t(t+2xW)W_t=0\ .$$
\end{cor}

\begin{rem} Partial derivations of $W$ and ${\overline W}$ are formal.
Proposition \ref{propanprop} (or a little work using 
assertion (ii) of Theorem \ref{thmeqW})
shows however that $W$ and ${\overline W}$
define analytic functions in suitable open subsets 
containing $(1,0)$ of $\mathbb C^*\times \mathbb C$.
\end{rem}

Section \ref{sectfewcompl} contains proofs for Theorem \ref{thmeqW},
\ref{thmcoeffid}, Corollary  \ref{cordiffeq} and a few
complements such as a homological interpretation for the coefficients
of the Laurent polynomials $w_n$ and features of the specialisations 
$W(-x,x)=1$ and $W(-x^{-1},x)=0$ of $W=W(t,x)$.

Section \ref{sectproofmain} contains the proof of
our main Theorem giving the following
closed formula for the coefficients of 
%the Laurent-polynomials 
$w_n(t)$.

\begin{thm} \label{thmwindcoeffs}
The Laurent-polynomials $w_n \in\mathbb N[t,t^{-1}]$ are
given by the following formulae: For every $n\geq 0$,
$$\begin{array}{l}
\displaystyle
w_{3n}=\sum_{k=0}^n{4n-2k\choose n+k}{2n+2k\choose 3k}
\frac{t^{n-2k}}{3k+1}\ ,\\
\displaystyle
w_{3n+1}=\sum_{k=0}^n{4n+2-2k\choose n+k}{2n+2k\choose 3k}
\frac{t^{n+1-2k}}{2n+1-k}\ ,\\
\displaystyle
w_{3n+2}=\sum_{k=0}^n{4n+2-2k\choose n+1+k}
{2n+2+2k\choose 3k+1}\frac{t^{n-2k}}{3k+2}\ .\\
\end{array}$$
\end{thm}

The outline of the proof for Theorem \ref{thmwindcoeffs}
is as follows: We display two partial differential equations 
which have at most a unique common solution.
We show then that the series defined by the formulae
of Theorem \ref{thmwindcoeffs} and the algebraic function
defined by assertion (ii) of Theorem \ref{thmeqW} 
are both common solutions of the differential equations mentionned above.
This implies that they coincide.

\begin{rem} Checkerboard triangulations of $P_n\supset e_*$
can be considered as the ${2n\choose n}\frac{1}{n+1}$ states of a 
spin model for the ``energy'' given by $t^{n_b-n_w}$. 
Theorem \ref{thmwindcoeffs} shows that this spin model is 
exactly solvable or integrable.
\end{rem}

\section{Main results for $d-$dissections}\label{sectmainddiss}

Dissections of polygons, also called cell-growth
problems, generalise triangulations coinciding 
with $2-$dissections.

Let $n\geq 0$ and $d\geq 2$ be two natural integers. We 
denote by $t_{d,n}$ the number of dissections of $P_{(d-1)n}$ 
into $n$ non-overlapping convex polygons such that every
polygon of the dissection has $(d+1)$ vertices which are all contained
in the set of vertices of $P_{(d-1)n}$. 
We call such a decomposition a {\it $d-$dissection}
of $P_{(d-1)n}$. The generating function 
$g_d=\sum_{n=0}^\infty t_{d,n}x^n\in\mathbb N[[x]]$ encoding the
numbers of all $d-$dissections satisfies the identity $g_d=1+xg_d^d$,
as shown by arguments similar to those used in the
proof of Theorem \ref{thmeqW}. 
The polynomial identity for $g_d$ follows 
also from the equivalence between $d-$dissections and
rooted plane $d-$regular trees. Section \ref{sectchir} 
explains this bijection in the case $d=2$. The
generalisation to arbitrary values of $d$ is straightforward.

The coefficients $t_{d,n}$ of $g_d=\sum_{n=0}^\infty t_{d,n}x^n$
are given by the following well-known 
result (see for example Formula 2.3 in \cite{HPR} or formula (5)
in \cite{Tak}).

\begin{thm} \label{thmcoeffgd}
The coefficients $t_{d,n}$ of
$g_d=1+x+dx^2+\dots=\sum_{n=0}^\infty t_{d,n}x^n$
defined by the equation $g_d=1+xg_d^d$ are given by the formula
$$t_{d,n}={dn\choose n}\frac{1}{(d-1)n+1}\ .$$
\end{thm}

For the sake of completeness, we include a short proof using 
Lagrange inversion at the beginning of Section \ref{sectgen}.

The choice of a marked edge 
$e_*\subset \partial P_{(d-1)n}$ turns a $d-$dissection 
$\tau$ of $P_{(d-1)n}$ into a {\it checkerboard $d-$dissection}
by partitioning all $n$ polygons involved in $\tau$ into two disjoint
subsets of $n_b$ black and $n_w$ white polygons such that the 
unique polygon $\Delta_*\in \tau$ containing the marked edge 
$e_*$ is black and such that edge-adjacent
polygons of $\tau$ have different colours. Associating 
the monomial $t^{n_b-n_w}$ to a checkerboard $d-$dissection
involving $n_b$ black and $n_w$ white polygons and summing the
monomials $t^{n_b-n_w}$ over all
$t_{d,n}$ checkerboard $d-$dissections of $P_{(d-1)n}\supset e_*$
yields again Laurent-polynomials $w_{d,n}\in\mathbb N[t,t^{-1}]$.

\begin{figure}[h]\label{fig1}
\epsfysize=1.4cm
\centerline{\epsfbox{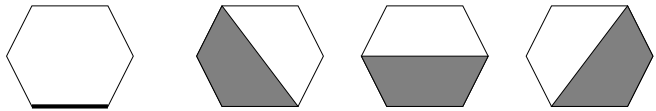}}
\caption{A hexagon with a marked edge and its $t_{3,2}=3$ 
contributions to $w_{3,2}=3$.}
\end{figure}

The proof of assertion (ii) in Theorem \ref{thmeqW}, modified
suitably, shows the following result.

\begin{thm} \label{thmgeqW} The generating series 
$W=\sum_{n=0}^\infty w_{d,n}x^n\in\mathbf R$
satisfies the algebraic equation
$$W=1+tx\left(1+t^{-1}xW^d\right)^d\ .$$
\end{thm}

Theorem \ref{thmcoeffid} and Corollary \ref{cordiffeq} have the following
easy and straightforward generalisations.

\begin{thm} \label{thmgcoeffid} We have the identity
$$((d-1)n+2+(d+1)j)(w_{d,n},t^{-j})=((d-1)n+2-(d+1)j)(w_{d,n},t^j)$$
where $(w_{d,n},t^j)$ denotes the coefficient of $t^j$ in 
$w_{d,n}=\sum_{j\in\mathbb Z} (w_{d,n},t^j)t^j\in\mathbb N[t,t^{-1}]$.
\end{thm}

Theorem \ref{thmgcoeffid} is of course equivalent to
$$(d-1)x{\overline W}_x+2{\overline W}+(d+1)t{\overline W}_t=(d-1)xW_x+2W-(d+1)tW_t$$
where $W=\sum_{n=0}^\infty w_{d,n}(t)x^n$ and ${\overline W}=
\sum_{n=0}^\infty w_{d,n}(t^{-1})x^n$.

\begin{cor} \label{corgdiffeq}
The generating series $W=\sum_{n=0}^\infty w_{d,n}x^n\in\mathbf R$
satisfies the partial differential equation
$$2t(1-W)+(d-1)x(dxW^{d-1}-t)W_x+(d+1)t(t+dxW^{d-1})W_t\ .$$
\end{cor}

Experimental observations suggest the following conjecture
for the coefficients of the Laurent-polynomials $w_{d,n}$.

\begin{conj}\label{conjcoeffwd} The Laurent-polynomials
$w_{d,n}\in \mathbb N[t,t^{-1}]$
are given by the following formulae: For every $n\geq 0$,
$$\begin{array}{l}
\displaystyle
w_{d,(d+1)n}=\sum_{k=0}^{(d-1)n}{d^2n-dk\choose n+k}{d(n+k)\choose (d+1)k}
\frac{t^{(d-1)n-2k}}{((d+1)k+1)}\\
\displaystyle
w_{d,(d+1)n+1}=\sum_{k=0}^{(d-1)n} {d^2n+d-dk\choose n+k}{d(n+k)\choose 
(d+1)k}\frac{t^{(d-1)n+1-2k}}{(dn+1-k)}
\end{array}$$
and 
$$\begin{array}{c}\displaystyle
w_{d,(d+1)n+j}=
\sum_{k=0}^{(d-1)n+j-2}{d^2n+d(j-1)-dk\choose n+1+k}
{dn+d+dk\choose (d+1)k+d-j+1}\ \cdot \\
\displaystyle \cdot \ \frac{t^{(d-1)n+j-2-2k}}{(d+1)k+d-j+2}
\end{array}$$
for $j$ such that $2\leq j\leq d$.
\end{conj}

Most steps in the proof of Theorem \ref{thmwindcoeffs} work
in the context of Conjecture \ref{conjcoeffwd}. Only a necessary
computation involving larger and larger polynomials for increasing
$d$ fails. Maple was however able to complete it for a few small values
of $d$ and we have:

\begin{thm} \label{thm234} Conjecture \ref{conjcoeffwd} holds
for $2\leq d\leq 6$.
\end{thm}

\section{A few complements and some easy proofs}\label{sectfewcompl}

\subsection{Proof of Theorem \ref{thmeqW}} A checkerboard 
triangulation $\tau$ of $P_n\supset e_*$ with $n\geq 1$ 
has a unique decomposition
into the marked black triangle $\Delta_*\in \tau$ containing $e_*$
and into a pair
$\tau_1,\tau_2$ of checkerboard triangulations with transposed colours 
of polygons $P_k,P_{n-1-k}$ obtained by cutting $\tau$ along 
$\Delta_*$ with a knife marking
all cuts. This construction is reversible and shows assertion (i).

Setting ${\overline W}=\sum_{n=0}^\infty w_n(t^{-1})x^n$, assertion (i) 
boils down to the equations 
$$W=1+tx{\overline W}^2\hbox{ and }{\overline W}=1+\frac{x}{t}W^2\ .$$
This shows that $W=1+tx(1+x/tW^2)^2$ and implies assertion (ii).
\hfill $\Box$

\subsection{Trees}\label{sectchir}

Triangulations of $P_n\supset e_*$ are well-known to be
in bijection with rooted plane $2-$regular trees, also 
called full binary trees, having $n$ interior vertices 
of outdegree $2$ and $n+1$ leaves, see for example the 
equivalence between the sets (i) and 
(vi) of Corollary 6.2.3 in \cite{St2}. 
The degenerate case $n=0$ corresponds to a rooted leaf
by convention. For $n\geq 1$,
the rooted tree $T_*(\tau)$ associated to a triangulation $\tau$
of $P_n\supset e_*$ has interior vertices given by the $n$ 
triangles of $\tau$, the root-vertex $v_*$ corresponding to the triangle 
$\Delta_*\supset e_*$,
and leaves given by the $n+1$ unmarked edges of $\partial P_n$.
Adjacency is given by pairs of vertices corresponding to
subsets in $\tau$ which intersect in a 
common edge. More precisely, two interior vertices of $T_*$ are 
adjacent if they correspond to edge-adjacent triangles of $\tau$.
An interior vertex $v$ of $T_*$ is adjacent to a leaf $l$ if
the triangle corresponding to $v$ contains the unmarked edge of 
$\partial P_n$ corresponding to $l$. The tree $T_*$ has $n_b$ interior
vertices at even distance from its root vertex $v_*$ and $n_w$ interior
vertices at odd distance from $v_*$. Otherwise stated, every 
tree is a bipartite graph and the numbers
$n_b,\ n_w$ count the interior vertices in the bipartite class
of the root vertex $v_*$ and in its complementary class.

\subsection{A homological interpretation for the coefficients of 
$w_n$}\label{subsechom}
 
Let $\tau$ be a checkerboard triangulation of $P_n\supset e_*$ where $n\geq 1$.
We denote as always by $\Delta_*\in\tau $ the unique black triangle containing
the marked edge $e_*\subset \partial P_{n}$. 
An easy induction on $n$
shows that there exists a unique continuous piecewise
affine map $\psi_\tau:P_n\longrightarrow \Delta_*$
which induces the identity-map on $\Delta_*$ and whose  
restriction to every triangle $\Delta\in\tau$ is an affine bijection.
We call $\psi_\tau$ the {\it folding map} since it maps, up to
piecewise affine transformations, the polygon $P_n$ onto $\Delta_*$
by folding $P_n$ along all $(n-1)$ interior edges of the
triangulation $\tau$ (where an edge $e\subset \Delta\in\tau$
is interior if $e\not\subset \partial P_n$). It is easy to check 
that $\psi_\tau$ preserves the orientation
on all black triangles and reverses the orientation on
all white triangles of $\tau$.

Orienting the boundaries $\partial P_n\subset P_n$ and 
$\partial \Delta_*\subset \Delta_*$ in the trigonometric counterclockwise
sense yields canonical isomorphisms between the first homology groups
$H_1(\partial P_n,\mathbb Z)\sim H_1(\partial \Delta_*,\mathbb Z)$
and the cyclic group $\mathbb Z$. The elements of 
$H_1(\partial P_n,\mathbb Z)$ or of $H_1(\partial \Delta_*,\mathbb Z)$
can be interpreted as ``winding coefficients'' of
closed loops contained in the boundary $\partial P_n$ or 
$\partial \Delta_*$ with respect to interior points
of the polygon $P_n$ or of the triangle $\Delta_*$.

The image $\psi_\tau(\partial \Delta)\subset \partial \Delta_*$
corresponds to the generator $1$ of $H_1(\partial \Delta_*,\mathbb Z)$
if $\Delta$ is a black triangle and to $-1$ if $\Delta$ is white.
This shows that the homomorphism
$$(\psi_\tau)_*:H_1(\partial P_n,\mathbb Z)\longrightarrow
H_1(\partial\Delta_*,\mathbb Z)$$
induced by the folding map is given by
$(\psi_\tau)_*(1)=n_b-n_w\in H_1(\partial \Delta_*,
\mathbb Z)$ where $n_b$ and $n_w$ are
the numbers of black and white triangles in
the checkerboard triangulation $\tau$.

Since every boundary edge $e\subset \partial P_n$ is contained in
a unique triangle (we assume $n\geq 1$)
of a checkerboard triangulation, it inherits 
a well-defined colour according to its inclusion
in the subset of black or white triangles.
In particular, the marked
edge $e_*$ of $\partial P_n$ is always black. The 
degenerate case $n=0$ corresponds by convention 
to a double edge consisting of a marked black edge and 
an unmarked white edge.

\begin{prop} \label{propformeb}
The numbers $e_b$ and $e_w$ of black and white edges in a 
checkerboard triangulation $\tau$ of $P_n\supset e_*$
are given by the formulae
$$\begin{array}{l}
\displaystyle e_b=\frac{n+2+3(n_b-n_w)}{2}=2n_b-n_w+1\\
\displaystyle e_w=\frac{n+2-3(n_b-n_w)}{2}=-n_b+2n_w+1\end{array}$$
with $n_b$ and $n_w$ denoting the number of black and white 
triangles in $\tau$.
\end{prop}

{\bf Proof} The homological interpretation of the coefficients in $w_n$
shows the identity $3(e_b-e_w)=n_b-n_w$. This, together with the trivial
equalities $e_b+e_w=n+2$ and $n_b+n_w=n$, implies the result.\hfill$\Box$

\subsection{The asymptotic mean-value of $n_b(\tau)-n_w(\tau)$}

The following result gives the 
asymptotic mean value $n_b-n_w$
of a uniform random triangulation for $P_n\supset e_*$.

\begin{prop} We have the equivalence
$$w'_n(1)\sim \frac{3}{8}{2n\choose n}\frac{1}{n+1}$$
for $n\rightarrow \infty$.

In particular, uniform random triangulations have the asymptotic
mean value $n_b-n_w=\frac{3}{8}$.
\end{prop}

{\bf Proof} 
Derivating the polynomial equation 
$$t(1+tx)-tW+2tx^2W^2+x^3W^4=0$$ 
for $W$ with respect to $t$, we get
$$1+2tx-W-tW_t+2x^2W^2+4tx^2WW_t+4x^3W^3W_t=0\ .$$
Eliminating $W$ yields the equation
$$\begin{array}{l}
\displaystyle
0=x(t+x)(8tx^2(2tx+1)-t+x)-t(8x^3(t^2+1)+12tx^2-t)W_t+\\
\displaystyle\qquad +2t^2x^2(64tx^3(t+x)+60tx^2-6t+48x^3+x)W_t^2+\\
\displaystyle\qquad +t^3x^3(256x^3(t^2+tx+1)+288tx^2-27t)W_t^4
\end{array}$$
for the partial derivation $W_t=\frac{\partial }{\partial t}W$.
The specialisation $t=1$ of the minimal polynomial for $W_t$ factorises
and yields the minimal polynomial 
$$x-\tilde W+x(4x+3)\tilde W^2=0$$
for
$$\tilde W=\sum_{n=0}^\infty w_n'(1)x^n=
\frac{1-(1+2x)\sqrt{(1-4x)}}{2x(4x+3)}=x+3x^3+4x^4+18x^5+\dots \ .$$
A straightforward computation using for example
$\lim_{n\rightarrow \infty} 
C_{n+1}/C_n=4$ for the Catalan numbers
$C_n={2n\choose n}\frac{1}{n+1}$ ends the proof.\hfill$\Box$

\subsection{Analytical properties of $W$}

\begin{prop} \label{propanprop} 
The series $W=W(t,x)$ is absolutely convergent
for $(t,x)\in\mathbb C^*\times \mathbb C$ such that 
$\max(\vert t^{1/3} x\vert,\ \vert t^{-1/3}x\vert) <\frac{1}{4}$.
\end{prop}

{\bf Proof} Non-negativity of the integers 
$e_b=\frac{n+2+3(n_b-n_w)}{2}$ and $e_w=\frac{n+2-3(n_b-n_w)}{2}$
given by Proposition \ref{propformeb} shows the bound 
$\vert n_b-n_w\vert\leq \frac{n+2}{3}$.
Positivity of all coefficients in $w_n$ yields the majorations
$$\vert w_n(t)\vert \leq w_n(1)M_n={2n\choose n}\frac{1}{n+1}M_n
\leq 4^n M_n$$
where $M_n=M_n(t)=\max(\vert t\vert^{(n+2)/3},
\vert t^{-1}\vert^{(n+2)/3})$. This implies the result easily.
\hfill$\Box$

\begin{rem} Analyticity of $W(t,x)$ in an open  neighbourhood of 
$\mathbb C^*\times\{0\}$ follows also from assertion (ii) in
Theorem \ref{thmeqW}. For fixed $t\in \mathbb C^*$,
the algebraic equation for $W(t,x)$ defines 
a unique analytic extension
$x\longmapsto W(t,x)$ of the evaluation $W(t,0)=1$
(the three remaining branches are singular at $x=0$).
\end{rem}

\begin{rem} Theorem \ref{thmwindcoeffs}. gives the exact bounds
on the degrees of monomials involved in $w_n$:
The Laurent polynomials
$w_{3n},w_{3n+2}$ contain no monomial of degree $>n$ or $<-n$ and
$w_{3n+1}$ contains no monomial of degree $>n+1$ or $<1-n$.
\end{rem}

\subsection{Proof of Theorem \ref{thmcoeffid} and
Corollary \ref{cordiffeq}}

A {\it direct automorphism} of a triangulation $\tau$ for $P_n$ is
a piecewise affine map inducing an orientation-preserving homeomorphism 
of $P_n$ which restricts to affine bijections between triangles of $\tau$.
The group of all such automorphisms is cyclic of order $\alpha^+(\tau)
\leq 3$.

A checkerboard triangulation $\tau$ of $P_n\supset e_*$ yields,
after unmarking the edge $e_*$, a contribution of 
$\frac{1}{\alpha^+(\tau)}e_b$ to the coefficient 
$t^{n_b-n_w}$ of $w_n$ and a contribution of 
$\frac{1}{\alpha^+(\tau)}e_w$ to the coefficient 
$t^{-n_b+n_w}$ of $w_n$. Theorem  \ref{thmcoeffid}. follows now easily
from Proposition \ref{propformeb}
for the numbers $e_b=\frac{n+2+3(n_b-n_w)}{2}$ and 
$e_w=\frac{n+2-3(n_b-n_w)}{2}$ of black and white edges in the 
checkerboard triangulation $\tau$ of $P_n\supset e_*$.
\hfill$\Box$

Corollary \ref{cordiffeq} follows from the equation
$x{\overline W}_x+2{\overline W}+3t{\overline W}_t=xW_x+2W-3tW_t$,
equivalent to Theorem \ref{thmcoeffid},
after elimination of ${\overline W},{\overline W}_t,
{\overline W}_x$ using the identity 
${\overline W}=1+\frac{x}{t}W^2$ occuring in the proof of 
Theorem \ref{thmeqW}
and its partial derivations
$W_t^-=-\frac{x}{t^2}W^2+2\frac{x}{t}WW_t$ and
$W_x^-=\frac{1}{t}W^2+2\frac{x}{t}WW_x$.

\subsection{Non-commutative edge polynomials} 
Labelling black and white edges in the boundary $\partial P_n$ 
of a checkerboard triangulation by $U$ and $V$, starting at
the marked edge $e_*$ and reading counterclockwise 
the labels of all $n+2$ edges
in $\partial P_n$, we get a word of length $n+2$, starting with $U$, 
in the alphabet $\{U,V\}$.
The sum over $n\in \mathbb N$ and over all checkerboard triangulations
of $P_n\supset e_*$ of these words defines thus a unique
non-commutative power-series $N\in \mathbb N\langle\!\langle
U,V\rangle\!\rangle$ in two free non-commuting variables.
The first few terms of $N$ are
$$\begin{array}{l}
\displaystyle
UV+U^3+UV^2U+U^2V^2+UV^4+UVU^3+U^2VU^2+U^3VU+U^4V+
\\
\displaystyle  \quad +2U^6+2(U^3V^3+U^2V^3U+UV^3U^2)+U^2VUV^2+UV^2UVU+\\
\displaystyle  \quad +U^2V^2UV+UV^2U^2V+UVU^2V^2+UVUV^2U+\dots\ .
\end{array}$$
Denoting by ${\overline N}=N(V,U)$ the series obtained by transposing the
variables $U,V$ of $N=N(U,V)$ we have the identities
$$\begin{array}{l}
N=UV+U\left(\frac{1}{V}{\overline N}\right)^2\\
{\overline N}=VU+V\left(\frac{1}{U}N\right)^2\end{array}$$
which imply
$$N=UV+U\left(U+\left(\frac{1}{U}N\right)^2\right)^2\ .$$

\subsection{The specialisations $t=\pm x$ and $t=\pm x^{-1}$}

\begin{prop} \label{propspec}
The specialisations $W(-x,x)$ and $W(-x^{-1},x)$
of $W(t,x)=\sum_{n=0}^\infty w_n(t)x^n$ are well defined and yield 
constant functions $W(-x,x)=1$ and $W(-x^{-1},x)=0$.
\end{prop}

This proposition amounts to annullation of all alternate row-sums 
except the first one and of all alternate column sums in the 
array $\mathcal A$ given by
$$\begin{array}{rrrrrrrrrrrrrrr}
1\\
1&2&1\\
&4&12&12&4\\
&2&30&100&140&90&22\\
&&28&280&980&1680&1540&728&140\\
&&9&360&2940&10584&20790&24024&16380&612&969\end{array}$$
and defined by writing the coefficients of $w_n$ suitably
along antidiagonals. 

Diagonal coefficients $1,2,12,100,980,\dots$
of $\mathcal A$ are constant terms in the Laurent 
polynomials $w_{2n}(t)$. Their generating series is the hypergeometric
function
$$\sum_{n=0}^\infty (w_{2n},t^0)x^n=\sum_{n=0}^\infty 
{2n\choose n}^2\frac{x^n}{n+1}$$
enumerating fair triangulations, see Theorem \ref{thmfair}
of the introduction.

The generating functions of the extremal sequences 
$$1,1,4,22,140,969,7084,\dots\hbox{ and }
1,2,9,52,340,2394,17710,\dots$$ of the array $\mathcal A$
are given by the algebraic hypergeometric functions
$$A=\sum_{n=0}^\infty {4n\choose n}\frac{x^n}{3n+1}=1+xB^2$$
and 
$$B=\sum_{n=0}^\infty {4n+2\choose n}\frac{x^n}{2n+1}=A^2,$$
cf. Section \ref{subsecEul} below.

{\bf Proof of Proposition \ref{propspec}} The upper bound
$\vert n_b-n_w\vert\leq \frac{n+2}{3}$
on the degrees of monomials with non-zero coefficient in $w_n$  
implies that only a finite number of non-zero coefficients of $W$
contribute to the monomial $x^m$ of $W(-x^{\pm 1},x)$.
They imply also $W(-x^{\pm 1},x)\in\mathbb Z[[x]]$.
The factorisations 
$$x(W-1)(1-x^2-x^2W+x^2W^2+x^2W^3)$$
corresponding to the specialisation $t=-x$ in assertion (ii) of
Theorem \ref{thmeqW} and
$$\frac{1}{x} W(1-2x^2W+x^4W^3)$$
corresponding to the specialisation $t=-x^{-1}$ in assertion (ii) of
Theorem \ref{thmeqW}, integrality of the ring 
$\mathbb Z[[x]]$ and a few easy verifications imply the result.
\hfill$\Box$ 

The generating functions $R(x)$ and $C(x)$ of row and column sums
of the array $\mathcal A$ considered above are given by
$$R(x)=W(\sqrt x,\sqrt x)=1+4x+32x^2+384x^3+5376x^4+82176x^5+\dots$$
satisfying the equation $1+x-R+2xR^2+xR^4=0$ and 
$$C(x)=W(\frac{1}{\sqrt x},\sqrt{x})=2+8x+80x^2+1024x^3+14848x^4+
231936x^5+\dots$$
satisfying $2-C+2xC^2+x^2C^4=0$.
The evaluations $W(-x,x)=1$ and $W(-x^{-1},x)=0$ show thus that 
row-sums or columns sums restricted to even elements of $\mathcal A$
(with indices starting at $0$) are given by the coefficients of
$\frac{1}{2}(R(x)+1)$ and $\frac{1}{2}C(x)$.

\begin{rem}
The equation $1+x-R+2xR^2+xR^4=0$ for $R=R(x)$ is equivalent to 
$x=\frac{R-1}{(1+R^2)^2}$. The formal power series
$R-1$ is thus the reciprocal series of the power series of the rational
function $y\longmapsto \frac{y}{(1+(1+y)^2)^2}$. 
Lagrange inversion (see for example Satz 2.4 
in \cite{Hen} or Theorem 5.4.2 in \cite{St2}) gives thus the formula
$$R=1+\sum_{n=1}^\infty \left(\sum_{k=0}^{\lfloor (n-1)/2\rfloor}
{2n\choose k}{2n-k\choose n-1-2k}\frac{1}{2^k}\right)\frac{(4x)^n}{n}\ .$$
\end{rem}

\subsection{Conjectural positivity properties}

The coefficient of $u^\alpha v^\beta$ involved in the
homogeneous polynomial $Q_n(u,v)=
\sqrt{uv}^n\ w_n\left(\sqrt{\frac{u}{v}}\right)
\in \mathbb N[u,v]$ of degree $n$ 
counts the number of checkerboard triangulations of  $P_n\supset e_*$ 
involving $\alpha$ black and $\beta$ white triangles.
Experimental observations suggest the following conjectural 
properties of the polynomials $Q_0,Q_1,\dots$.

\begin{conj} (i)
The polynomials $Q_0(u,1),Q_1(u,1),\dots\in \mathbb N[u]$
have only real roots (which are $\leq 0$ since all coefficients
are positive) and the non-zero roots of $Q_n(u,1)$
interleave the non-zero roots of $Q_{n+3}(u,1)$.

\ \ (ii) The symmetric polynomial $Q_n(u,v)+Q_n(v,u)$
is of the form $R_n(u+v,uv)$ where $R_n(e_1,e_2)\in\mathbb Z[e_1,e_2]$
involves only positive coefficients.
\end{conj}

\subsection{Eulerian triangulations}\label{subsecEul}

Theorem \ref{thmwindcoeffs} or fairly elementary
generating series manipulations show that $P_{3n+1}\supset e_*$
has exactly ${4n+2\choose n}\frac{1}{2n+1}$ checkerboard triangulations
such that the boundary $\partial P_{3n+1}$ of $P_{3n+1}$
is contained in the subset of all black triangles. 
Such a triangulation $\tau$ could be called ``Eulerian'' since 
all vertices of $P_{3n+1}$ have even degree in the planar graph 
defined by $\tau$.
Equivalently, a triangulation $\tau$ of $P_n$ is Eulerian 
if every vertex of $P_n$ belongs to an odd number of triangles
in $\tau$.

\begin{rem}
The boundary $\partial P_{(d+1)n+1}$ of $P_{(d+1)n+1}$ seems to be
contained in the subset of black polygons for exactly
${d^2n+d\choose n}\frac{1}{dn+1}$ checkerboard $d-$dissections.
Such dissections define again Eulerian graphs and have the equivalent
property that every vertex of $P_{(d+1)n+1}$ is contained 
in an odd number of dissecting polygons. 
\end{rem}

\section{Proof of Theorem \ref{thmwindcoeffs}}\label{sectproofmain}

The idea for proving Theorem \ref{thmwindcoeffs} is as follows:

We exhibit two partial differential operators $D=D_L-D_R$ and 
$\tilde D=\tilde D_L-\tilde D_R$ such that the 
two associated partial differential equations $DF=0$ 
and $\tilde DF=0$ have at most a unique common formal solution 
satisfying the initial condition $F\equiv 1+tx\pmod{x^2}$ in
the algebra $\mathbf R=\left(\mathbb Q[t,t^{-1}]\right)[[x]]$
of formal power series in $x$ with coefficients in $\mathbb Q[t,t^{-1}]$.

We show then that we have $D\tilde W=\tilde D\tilde W=0$
for the formal series $\tilde W\equiv 1+tx\pmod{x^2}$
defined by the formulae of Theorem \ref{thmwindcoeffs}.

Finally, we consider the algebraic equation 
$$t(1+tx)-tW+2tx^2W^2+x^3W^4=0$$
for the generating series $W=\sum_{n=0}^\infty w_nx^n$
enumerating checkerboard triangulations. We show that this equation
has a unique solution in $\mathbf R$ which is thus 
given by $W$
satisfying the initial condition $W\equiv 1+tx\pmod{x^2}$.
We show then that $W$
satisfies the equations $DW=\tilde DW=0$. This implies $W=\tilde W$
and establishes Theorem \ref{thmwindcoeffs}.

Consider the four linear partial differential operators
$$\begin{array}{lcl}
\displaystyle D_L&\displaystyle =&
\displaystyle 4tx\left(1+t\frac{\partial}{\partial t}+x\frac{\partial}
{\partial x}\right)\left(2-3t\frac{\partial}{\partial t}+x\frac{\partial}
{\partial x}\right)\\
&\displaystyle =&\displaystyle 4tx\left(2-4t\frac{\partial}{\partial t}
+4x\frac{\partial}{\partial x}-3t^2\frac{\partial^2}{\partial 
t^2}-2tx\frac{\partial^2}{\partial t\partial x}+x^2
\frac{\partial^2}{\partial x^2}\right)\\
\displaystyle D_R&\displaystyle =&
\displaystyle \left(3t\frac{\partial}{\partial t}+x\frac{\partial}
{\partial x}\right)\left(-2+3t\frac{\partial}{\partial t}+x\frac{\partial}
{\partial x}\right)\\
&\displaystyle =&\displaystyle 3t\frac{\partial}{\partial t}-x\frac{\partial}{\partial x}+
9t^2\frac{\partial^2}{\partial t^2}+6tx\frac{\partial}{\partial t\partial x}
+x^2\frac{\partial^2}{\partial x^2}\\
\displaystyle \tilde D_L&\displaystyle =&\displaystyle 
4x\left(1-t\frac{\partial}{\partial t}+x\frac{\partial}
{\partial x}\right)\left(3t\frac{\partial}{\partial t}+x\frac{\partial}
{\partial x}\right)\\
&\displaystyle =&\displaystyle 4x\left(
2x\frac{\partial}{\partial x}-3t^2\frac{\partial^2}{\partial t^2}
+2tx\frac{\partial^2}{\partial t\partial x}+x^2\frac{\partial^2}
{\partial x^2}\right)\\
\displaystyle \tilde D_R&\displaystyle =&
\displaystyle t\left(2-3t\frac{\partial}{\partial t}+x\frac{\partial}
{\partial x}\right)\left(-3t\frac{\partial}{\partial t}+x\frac{\partial}
{\partial x}\right)\\
&\displaystyle =&\displaystyle t\left(3t\frac{\partial}{\partial t}
+3x\frac{\partial}{\partial x}+9t^2\frac{\partial^2}{\partial t^2}
-6tx\frac{\partial^2}{\partial t\partial x}+x^2\frac{\partial^2}
{\partial x^2}\right)
\end{array}$$

\begin{prop} \label{propsoluni} The two partial differential equations 
$$D_LF=D_RF$$
and 
$$\tilde D_LF=\tilde D_RF$$
have at most a unique common solution $F\in\mathbf R=(\mathbb 
Q[t,t^{-1}])[[x]]$
which satisfies the initial condition $F\equiv 1+tx\pmod {x^2}$.
\end{prop}

{\bf Proof} The formulae 
$$\begin{array}{l}
\displaystyle D_L(t^jx^m)=4(1+j+m)(2-3j+m)t^{j+1}x^{m+1}\\
\displaystyle D_R(t^jx^m)=(3j+m)(-2+3j+m)t^jx^m\\
\displaystyle \tilde D_L(t^jx^m)=4(1-j+m)(3j+m)t^jx^{m+1}\\
\displaystyle \tilde D_R(t^jx^m)=(2-3j+m)(-3j+m)t^{j+1}x^m
\end{array}$$
show that a coefficient $(F,t^jx^m)$ of a common solution $F$
is determined by the coefficients
$(F,t^{j-1}x^{m-1})$ and $(F,t^{j+1}x^{m-1})$
except if $(j,m)\in\mathbb Z\times \mathbb N$  is among the four common roots
$(0,0),(\frac{1}{3},-1),(\frac{1}{3},1),(\frac{2}{3},0)$
of the two polynomials
$$\begin{array}{l}
\displaystyle (3j+m)(-2+3j+m)\ ,\\
\displaystyle (2-3j+m)(-3j+m)\ .\end{array}$$
This shows $(j,m)=(0,0)$ and implies that 
such a solution $F$ is either non-existent or uniquely defined by 
the initial condition $F\equiv 1+tx\pmod{x^2}$.
\hfill$\Box$

\begin{prop} \label{proptWsols}
We have 
$$D_L\tilde W=D_R\tilde W$$ and
$$\tilde D_L\tilde W=\tilde D_R\tilde W$$ 
for the series $\tilde W=\sum_{n=0}^\infty \tilde w_nx^n\equiv 1+tx
\pmod{x^2}$
defined by the formulae
$$\begin{array}{l}
\displaystyle
\tilde w_{3n}=\sum_{k=0}^n{4n-2k\choose n+k}{2n+2k\choose 3k}
\frac{t^{n-2k}}{3k+1}\ ,\\
\displaystyle
\tilde w_{3n+1}=\sum_{k=0}^n{4n+2-2k\choose n+k}{2n+2k\choose 3k}
\frac{t^{n+1-2k}}{2n+1-k}\ ,\\
\displaystyle
\tilde w_{3n+2}=\sum_{k=0}^n{4n+2-2k\choose n+1+k}
{2n+2+2k\choose 3k+1}\frac{t^{n-2k}}{3k+2}\ .\\
\end{array}$$
\end{prop}

\begin{lem} \label{lemDLDR} (i) We have
$$\begin{array}{l}
D_L(x^{3n}\tilde w_{3n}(t))=8x^{3n+1}\sum_{k=0}^n\frac{(4n+1-2k)!\ (2n+2k)!
\ t^{n+1-2k}}
{(n+k)!\ (3n-3k)!\ (3k)!\ (2n-k)!}\\
D_L(x^{3n+1}\tilde w_{3n+1}(t))=\\
\qquad 8x^{3n+2}\sum_{k=1}^n
\frac{(4n+1-2(k-1))!\ (2n+2+2(k-1))!\ t^{n-2(k-1)}}{(n+1+(k-1))!\ 
(3n-1-3(k-1))!\ (3(k-1)+2)!\ (2n-(k-1))!}\\
D_L(x^{3n+2}\tilde w_{3n+2}(t))=\\
\qquad 4x^{3(n+1)}\sum_{k=0}^n
\frac{(4(n+1)-2k)!\ (2(n+1)+2k)!\ t^{n+1-2k}}{(n+1+k)!\ (3(n+1)-3k-2)!\
(3k+1)!\ (2(n+1)-k)!}
\end{array}$$
and
$$\begin{array}{l}
D_R(x^{3n}\tilde w_{3n})=4x^{3n}\sum_{k=0}^{n-1}
\frac{(4n-2k)!\ (2n+2k)!\ t^{n-2k}}{(n+k)!\ (3n-3k-2)!\ (3k+1)!\ 
(2n-k)!}\\
D_R(x^{3n+1}\tilde w_{3n+1}(t))=8x^{3n+1}\sum_{k=0}^n
\frac{(4n+1-2k)!\ (2n+2k)!\ t^{n+1-2k}} {(n+k)!\ (3n-3k)!\ (3k)!\ (2n-k)!}\\
D_R(x^{3n+2}\tilde w_{3n+2}(t))=8x^{3n+2}\sum_{k=0}^{n-1}
\frac{(4n+1-2k)!\ (2n+2+2k)!\ t^{n-2k}}{(n+1+k)!\ (3n-1-3k)!\ (3k+2)!\ 
(2n-k)!}\\
\end{array}$$

\ \ (ii) We have
$$\begin{array}{l}
\tilde D_L(x^{3n}\tilde w_{3n}(t))=8x^{3n+1}\sum_{k=0}^{n-1}\frac{(4n-2k)!\ (2n+1+2k)!
\ t^{n-2k}}
{(n+k)!\ (3n-3k-1)!\ (3k+1)!\ (2n-k)!}\\
\tilde D_L(x^{3n+1}\tilde w_{3n+1}(t))=8x^{3n+2}\sum_{k=0}^n
\frac{(4n+2-2k)!\ (2n+1+2k)!\ t^{n+1-2k}}{(n+k)!\ (3n+1-3k)!\ (3k)!\ (2n+1-k)!}
\\
\tilde D_L(x^{3n+2}\tilde w_{3n+2}(t))=\\
\qquad 8x^{3(n+1)}\sum_{k=0}^{(n+1)-1}
\frac{(4(n+1)-2-2k)!\ (2(n+1)+1+2k)!\ t^{(n+1)-1-2k}}{((n+1)+k)!
\ (3(n+1)-3k-3)!\ (3k+2)!\ (2(n+1)-1-k)!}
\end{array}$$
and
$$\begin{array}{l}
\tilde D_R(x^{3n}\tilde w_{3n})=8x^{3n}\sum_{k=1}^{n}
\frac{(4n-2-2(k-1))!\ (2n+1+2(k-1))!\ t^{n-1-2(k-1)}}{(n+(k-1))!\ 
(3n-3(k-1)-3)!\ (3(k-1)+2)!\ (2n-1-(k-1))!}\\
\tilde D_R(x^{3n+1}\tilde w_{3n+1}(t))=\\
\qquad 8x^{3n+1}\sum_{k=1}^n
\frac{(4n-2(k-1))!\ (2n+1+2(k-1))!\ t^{n-2(k-1)}} {(n+(k-1))!\ (3n-3(k-1)-1)!\ (3(k-1)+1)!\ (2n-(k-1))!}\\
\tilde D_R(x^{3n+2}\tilde w_{3n+2}(t))=8x^{3n+2}\sum_{k=0}^{n}
\frac{(4n+2-2k)!\ (2n+1+2k)!\ t^{n+1-2k}}{(n+k)!\ (3n+1-3k)!\ (3k)!\ 
(2n+1-k)!}\\
\end{array}$$
\end{lem}

{\bf Proof} Elementary and tedious verifications left to the reader.
\hfill$\Box$

{\bf Proof of Proposition \ref{proptWsols}}
Follows easily from Lemma \ref{lemDLDR}.\hfill$\Box$

%\subsection{A polynomial identity}

Consider the ideal $\mathcal I=
\mathcal (P,P_t,P_x,P_{tt},P_{tx}P_{xx})$ generated by
$$\begin{array}{l}
P=t(1+tx)-tF+2tx^2F^2+x^3F^4,\\
P_t=1+2tx-F-tF_t+2x^2F^2+4tx^2FF_t+4x^3F^3F_t,\\
P_x=t^2-tF_x+4txF^2+4tx^2FF_x+3x^2F^4+4x^3F^3F_x,\\
P_{tt}=2x-2F_t-tF_{tt}+8x^2FF_t+4tx^2F_t^2+4tx^2FF_{tt}+\\
\qquad +12x^3F^2F_t^2+4x^3F^3F_{tt},\\
P_{tx}=2t-F_x-tF_{tx}+4xF^2+4x^2FF_x+8txFF_t+4tx^2F_tF_x+\\
\quad +4tx^2FF_{tx}+12x^2F^3F_t+12x^3F^2F_tF_x+4x^3F^3F_{tx},\\
P_{xx}=-tF_{xx}+4tF^2+16txFF_x+4tx^2F_x^2+4tx^2FF_{xx}+6xF^4+\\
\quad +24x^2F^3F_x+12x^3F^2F_x^2+4x^3F^3F_{xx}\end{array}$$
of the free polynomial algebra $\mathbb Q[t,x,F,F_t,F_x,F_{tt},F_{tx},F_{xx}]$
in eight variables $t,x,F,F_t,F_x,F_{tt},F_{tx},F_{xx}$.

We introduce moreover the polynomials
$$\begin{array}{l}
\displaystyle Q=4tx\left(2F-4tF_t+4xF_x-3t^2F_{tt}-2txF_{tx}
+x^2F_{xx}\right)\\
\displaystyle \qquad -(3tF_t-xF_x+9t^2F_{tt}+6txF_{tx}+x^2F_{xx})\\
\tilde Q=4x\left(2xF_x-3t^2F_{tt}+2txF_{tx}+x^2F_{xx}\right)\\
\displaystyle 
-t\left(3tF_t+3xF_x+9t^2F_{tt}-6txF_{tx}+x^2F_{xx}\right)\end{array}$$
and
$$K=\frac{\partial P}{\partial F}=-t+4tx^2F+4x^3F^3\ .$$

\begin{lem}\label{lemK3Q} We have $K^3Q\in\mathcal I$
and $K^3\tilde Q\in\mathcal I$.
\end{lem}

{\bf Proof} 
Set
$$Q_1=KQ-((Q,F_{tt})P_{tt}+(Q,F_{tx})P_{tx}+(Q,F_{xx})P_{xx})
\equiv KQ\pmod{\mathcal I}$$
and 
$$\tilde Q_1=K\tilde Q-((\tilde 
Q,F_{tt})P_{tt}+(\tilde Q,F_{tx})P_{tx}+(\tilde Q,F_{xx})P_{xx})
\equiv K\tilde Q\pmod{\mathcal I}$$
where $(R,F_{tt}),(R,F_{tx}),(R,F_{xx})\in\mathbb Q[t,x,F,F_t,F_x]$
are the coefficients of $F_{tt},F_{tx},F_{xx}$ of $R=Q$ or $R=\tilde Q$.
Since the three 
polynomials $KF_{tt}-P_{tt},KF_{tx}-P_{tx},KF_{xx}-P_{xx}$
are elements of $\mathbb Q[t,x,F,F_t,F_x]$ and 
since $Q$ and $\tilde Q$ are of degree $1$ with respect to the variables
$F_{tt},F_{tx},F_{xx}$ we have the inclusions
$Q_1,\tilde Q_1\in\mathbb Q[t,x,F,F_t,F_x]$.

Similarly, we have $KF_t-P_t,KF_x-P_x\in\mathbb Q[t,x,F]$.
Since $Q_1,\tilde Q_1$ are of degree $2$ with respect to the variables
$F_t,F_x$, substituting $F_t$ by $F_t-\frac{1}{K}P_t$ and $F_x$
by $F_x-\frac{1}{K}P_x$ in $K^2Q_1,K^2\tilde Q_1$ yields
polynomials $Q_0\equiv K^2Q_1\equiv K^3Q\pmod{\mathcal I}$
and $\tilde Q_0\equiv K^2\tilde Q_1\equiv K^3\tilde Q$ 
such that $Q_0,\tilde Q_0\in\mathbb Q[t,x,F]$.
A computation (using for example Maple)
shows that $P$ divides $Q_0$ and $\tilde Q_0$ which ends the proof.
\hfill $\Box$.

\begin{prop} \label{propWunique} The algebraic equation
$$t(1+tx)-tW+2tx^2W^2+x^3W^4=0$$
has a unique solution $W$ in $\mathbf R=(\mathbb Q[t,t^{-1}])[[x]]$.

This solution is a common solution of the two partial differential
equations 
$$\begin{array}{l}
\displaystyle
D_L W=D_R W\\
\displaystyle\tilde D_L W=\tilde D_R W\end{array}$$ 
and satisfies the initial condition $W=1+tx\pmod{x^2}$.
\end{prop}

{\bf Proof of Theorem \ref{thmwindcoeffs}}
Follows from Propositions \ref{propsoluni},
\ref{proptWsols} and \ref{propWunique}.\hfill$\Box$

{\bf Proof of Proposition \ref{propWunique}} The solution of
$$t(1+tx)-tW+2tx^2W^2+x^3W^4=0$$
is the unique fixpoint in $\mathbf R=\left(\mathbb C[t,t^{-1}]
\right)[[x]]$
of the attracting map
$Z\longmapsto 1+tx(1+x/tZ^2)^2$.

Since the algebra $\mathbf R$ containing the solution $W$
considered above is  a differential algebra for both partial derivations
$\frac{\partial}{\partial t}$ and $\frac{\partial}{\partial x}$, 
we can consider the homomorphism of algebras
$$\varphi:
\mathbb Q[t,x,F,F_t,F_x,F_{tt},F_{tx},F_{xx}]\longrightarrow \mathbf R$$
defined by 
$$\varphi(t)=t,\ \varphi(x)=x,\ \varphi(F_{t^\alpha x^\beta})=
\frac{\partial^{\vert \alpha+\beta\vert}}{\partial t^\alpha 
\partial x^\beta}W,\ \alpha+\beta\leq 2\ .$$
Since the polynomials $P_{t^\alpha x^\beta}$
are formally given by $\frac{\partial^{\alpha+\beta}}{\partial t^\alpha
\partial x^\beta}P$, we have
$\varphi(P_{t^\alpha x^\beta})=0$
for all $\alpha,\beta\in\mathbb N$. This implies the inclusion $\mathcal 
I\subset \mathop{ker}(\varphi)$ for the ideal $\mathcal I=(
P,P_t,P_x,P_{tt},P_{tx},P_{xx})$.

The identities
$\varphi(Q)=(D_L-D_R)W$ and $\varphi(\tilde Q)=(\tilde D_L-\tilde D_R)W$
and Lemma \ref{lemK3Q} imply thus the equalities 
$\varphi(K)^3(D_L-D_R)W=0$ and $\varphi(K)^3(\tilde D_L-\tilde D_R)W=0$.
Since the algebra $\mathbf R$ has no zero divisors,
and since $\varphi(K)\equiv -t\pmod{x^2}$, we have
$(D_L-D_R)W=0$ and $(\tilde D_L-\tilde D_R)W=0$ with 
$W\equiv 1+tx\pmod{x^2}$.\hfill$\Box$

\begin{rem} Lemma \ref{lemK3Q} can be replaced by the inclusions
$Q^2,\tilde Q\in\mathcal I$ which can be checked by computing 
a Gr\"obner basis for $\mathcal I$. I thank Bernard Parisse who 
did the necessary computations using Xcas and CoCoa. 
The computation of a Gr\"obner 
basis is however typically quite long (several minutes
in the above case) while the computations used for the previous proof of 
Lemma \ref{lemK3Q} are immediate.
\end{rem}

\begin{rem} Corollary \ref{cordiffeq} results also from the 
following computation
which is analogous to the proof of Lemma \ref{lemK3Q}: Replacing
$F_t$ by $F_t-\frac{1}{K}P_t$ and 
$F_x$ by $F_x-\frac{1}{K}P_x$ in
$$K\big(2t(1-F)+x(2xF-t)F_x+3t(t+2xF)F_t\big)$$
(with $K=-t+4tx^2F+4x^3F^3$ as above) we get
$-(5t+6xF)P$
which implies Corollary \ref{cordiffeq}.
\end{rem}

\section{(Partial) proofs for $d-$dissections}\label{sectgen}

{\bf Proof of Theorem \ref{thmcoeffgd}}
Setting $G=G(x)=xg_d(x^{d-1})$, the computation
$$G-G^d=x\left(g_d(x^{d-1})-x^{d-1}g_d(x^{d-1})^{d-1}\right)=x$$
shows that $G(x)=x+\dots$ is the reciprocal function of 
$y\longmapsto y-y^d$. Since the coefficient $\gamma_n$ of 
$g_d(x)=\sum_{n=0}^\infty \gamma_n x^n$
equals the coefficient
of $x^{n(d-1)+1}$ in $G(x)=xg_d(x^{d-1})$, Lagrange inversion 
(which states that $n\sigma_n\in\mathbb C$ is given by 
the coefficient of $y^{n-1}$
in $(y/r(y))^n\in\mathbb C[[y]]$ if $\sum_{n=1}^\infty \sigma_n (r(y))^n=y$
for $r(y)=\sum_{n=1}^\infty \rho_ny^n\in\mathbb C[[y]]$, see eg. Satz 2.4 
in \cite{Hen} or Theorem 5.4.2 in \cite{St2}),
shows that $(n(d-1)+1)\gamma_n$ is given by the coefficient of
$y^{n(d-1)}$ in
$$\begin{array}{rcl}
\displaystyle \left(\frac{y}{y-y^d}\right)^{n(d-1)+1}&
\displaystyle =&\displaystyle\sum_{j=0}^\infty
{-n(d-1)-1\choose j} (-1)^j y^{j(d-1)}\\
&\displaystyle =&\displaystyle \sum_{j=0}^\infty
{n(d-1)+j\choose j}y^{j(d-1)}\end{array}$$
and we have thus
$$\gamma_n={nd\choose n}\frac{1}{n(d-1)+1}$$
as required.\hfill$\Box$

{\bf Proof of Theorem \ref{thmgeqW}} Introducing $\tilde W=\tilde W(t,x)=
W(t^{-1},x)\in\mathbf R$ and decomposing for $n\geq 1$ 
a $d-$dissection $\tau$ of $D_{(d-1)n}\supset
e_*$
along the distinguished polygon $\Delta_*\in\tau$ containing $e_*$, 
we get the equations
$$\begin{array}{l}
\displaystyle W_d=1+tx\tilde W_d^d\\
\displaystyle \tilde W_d=1+t^{-1}xW_d^d\ .\end{array}$$
Eliminating $\tilde W$ yields the algebraic equation 
$$W_d=1+tx\left(1+t^{-1}xW_d^d\right)^d$$
of Theorem \ref{thmgeqW} for $W$.\hfill$\Box$

\subsection{A partial ``proof'' for Conjecture \ref{conjcoeffwd}}

All steps but one in the proof of Theorem \ref{thmwindcoeffs} 
work for arbitrary $d\geq 2$ and yield thus almost a proof 
of Conjecture \ref{conjcoeffwd}.
Failure occurs in the necessary machine computations which get more
and more complicated for increasing values of $d$. 
We were however able to complete them for a few small values of $d$
and thus to establish Theorem \ref{thm234}.

Given a constant $d\geq 2$, we consider the four following
partial differential operators
$$D_L=dtx\left(1-\frac{(d+1)}{2}t\frac{\partial}{\partial t}+
\frac{(d-1)}{2}x\frac{\partial}{\partial x}\right)
\ \prod_{j=1}^{d-1}\left(j+\frac{d}{2}t\frac{\partial}{\partial t}+
\frac{d}{2}x\frac{\partial}{\partial x}\right)$$

$$D_R=\prod_{j=0}^{d-1}\left(-j+\frac{(d+1)}{2}t\frac{\partial }{\partial t}
+\frac{(d-1)}{2}x\frac{\partial}{\partial x}\right)$$

$$\tilde D_L=dx\left(\frac{(d+1)}{2}t\frac{\partial}{\partial t}+
\frac{(d-1)}{2}x\frac{\partial}{\partial x}\right)
\ \prod_{j=1}^{d-1}\left(j-\frac{d}{2}t\frac{\partial}{\partial t}+
\frac{d}{2}x\frac{\partial}{\partial x}\right)$$

$$\tilde D_R=t\prod_{j=0}^{d-1}\left(1-j-\frac{(d+1)}{2}t\frac{\partial }{\partial t}
+\frac{(d-1)}{2}x\frac{\partial}{\partial x}\right)$$

\begin{prop} The two partial differential equations 
$$D_LF=D_RF$$
and 
$$\tilde D_LF=\tilde D_RF$$
defined by the previous partial differential operators 
have at most a unique common solution $F\in\mathbf R$
satisfying the initial condition $F\equiv 1+tx\pmod {x^2}$.
\end{prop}

{\bf Proof} As in the proof of the special case 
$d=2$ (see Proposition \ref{propsoluni}),
a coefficient $(F,t^jx^m)$ of a common solution $F$
is determined by the coefficients
$(F,t^{j\pm1}x^{m-1})$
except if $D_R(t^jx^m)=\tilde D_R(t^jx^m)=0$.
Such a pair of integers $(j,m)$ satisfies the
two linear equations
$$-a+\frac{d+1}{2}j+\frac{d-1}{2}m=1-b-\frac{d+1}{2}j+\frac{d-1}{2}m=0$$
for some $a,b\in\{0,\dots,d-1\}$.
Adding these two linear equations we have
$$(d-1)m=a+b-1\leq 2d-3$$
which shows $m\in\{0,1\}$. 
A coefficient of $t^j$ or of $t^jx$ in a solution $F$ is however
prescribed by the initial condition $F\equiv 1+tx\pmod{x^2}$.
\hfill $\Box$

\begin{prop} \label{propgtWsols}
We have 
$$D_L\tilde W=D_R\tilde W$$ and
$$\tilde D_L\tilde W=\tilde D_R\tilde W$$ 
for the series $\tilde W=\sum_{n=0}^\infty \tilde w_{d,n}x^n\equiv 1+tx
\pmod{x^2}$
defined by the formulae given in Conjecture \ref{conjcoeffwd}.
\end{prop}

{\bf Proof} Follows from the formulae

$$\begin{array}{l}
D_L(x^{(d+1)n+1}\tilde w_{(d+1)n+1})=dx^{(d+1)n+2}\sum_{k=1}^{(d-1)n}\\
\quad \frac{(d^2n+d-d(k-1)-1)!\  (dn+d+d(k-1))!\ t^{(d-1)n-2(k-1)}  }
{(n+(k-1)+1)!\
((d^2-1)n-(d+1)(k-1)-1)!\ (d(k-1)+d+(k-1))!
\ (dn-(k-1))!}\end{array}$$

$$\begin{array}{l}
D_L(x^{(d+1)n+j}\tilde w_{(d+1)n+j})=dx^{(d+1)n+(j+1)}\sum_{k=0}^{(d-1)n+j-2}\\
\quad \frac{(d^2n+dj-dk-1)!\  (dn+d(k+1))!\ 
t^{(d-1)n+(j+1)-2(k+1)}  }{(n+(k+1))!\
((d^2-1)n+dj-(d+1)(k+1))!\ ((d+1)(k+1)-j)!\ (dn+(j+1)-(k+1)-1)!}
\end{array}$$

$$\begin{array}{l}
D_R(x^{(d+1)n+1}\tilde w_{(d+1)n+1})=x^{(d+1)n+1}\sum_{k=0}^{(d-1)n}\\
\quad \frac{(d^2n+d-dk)!\ (dn+dk)!\ t^{(d-1)n+1-2k}}{(n+k)!\
((d^2-1)n-(d+1)k)!\ ((d+1)k)!\ (dn+1-k)!}\end{array}$$

$$\begin{array}{l}
D_R(x^{(d+1)n+j}\tilde w_{(d+1)n+j})=x^{(d+1)n+j}\sum_{k=0}^{(d-1)n+j-3}\\
\quad \frac{(d^2n+d(j-1)-dk)! \  (dn+d(k+1))!\ t^{(d-1)n+j-2-2k} }{(n+k+1)!\
((d^2-1)n+d(j-2)-(d+1)k-1)!\
((d+1)k+d-j+2)!\ (dn+j-k-1)!}\end{array}$$

$$\begin{array}{l}
\tilde D_L(x^{(d+1)n+1}\tilde w_{(d+1)n+1})=dx^{(d+1)n+2}\sum_{k=0}^{(d-1)n}\\
\quad \frac{(d^2n+d-dk)!\  (dn+d+dk-1)!\ t^{(d-1)n+1-2k}   }{(n+k)!\
((d^2-1)n+d-(d+1)k-1)!\ (d+1)k)!\ (dn+1-k)!}\end{array}$$

$$\begin{array}{l}\tilde D_L(x^{(d+1)n+j}\tilde w_{(d+1)n+j})=
dx^{(d+1)n+(j+1)}\sum_{k=0}^{\omega_d(j)}\\
\quad\frac{(d^2n+d(j-1)-dk)!\   (dn+2d+dk-1))!\ t^{(d-1)n+j-2-2k}  }{(n+k+1)!\
((d^2-1)n+dj-(d+1)(k+1)-1)!\ 
((d+1)k+d-j+2)!\ (dn+j-k-1)!}\end{array}$$
where $\omega_d(j)=(d-1)n+j-2$ for $j\in\{2,3,\dots,d\}$ and $\omega_d(d+1)=
(d-1)n+d-2$,

$$\begin{array}{l}
\tilde D_R(x^{(d+1)n+1}\tilde w_{(d+1)n+1})=x^{(d+1)n+1}\sum_{k=1}^{(d-1)n}\\
\quad \frac{(d^2n+d-dk)!\ (dn+dk)!\ t^{(d-1)n+2-2k}}{(n+k)!\
((d^2-1)n+d-(d+1)k)!\ ((d+1)k-d)!\ (dn+1-k)!}\end{array}$$

$$\begin{array}{l}
\tilde D_R(x^{(d+1)n+j}\tilde w_{(d+1)n+j})=x^{(d+1)n+j}
\sum_{k=\alpha_d(j)}^{(d-1)n+j-2}\\
\quad \frac{(d^2n+d(j-1)-dk)!\  (dn+d+dk)!\ t^{(d-1)n+j-1-2k}  }{(n+k+1)!\
((d^2-1)n+d(j-1)-(d+1)k-1)!\ ((d+1)k-j+2)!\ (dn+j-k-1)!}\end{array}$$
where $\alpha_d(2)=0$ and $\alpha_d(j)=1$ for $j\in\{3,\dots,d+1\}$,
corresponding to Lemma \ref{lemDLDR} and from the observation 
that the formula of Conjecture \ref{conjcoeffwd}
for $\tilde w_{d,(d+1)n+j}$ with $j=d+1$ coincides with
the formula for $\tilde w_{d,(d+1)(n+1)}$.\hfill $\Box$

\subsection{An algebraic reformulation}

The integer $d\geq 2$ is again fixed in this Section.
We denote by $\mathcal A_h=\mathbb Q[t,x,\left(F_{t^\alpha x^\beta}\right)_{
\alpha+\beta\leq h}]$
the free algebra generated by $t,x$ and all partial derivations 
$F_{t^\alpha x^\beta}=\frac{\partial^{\alpha+\beta}}{\partial t^\alpha\
\partial x^\beta}F$ of order $\alpha+\beta\leq h$ of an unknown 
analytic function $F=F(t,x)$. We suppose
that there are no algebraic relations among partial derivations of $F$.

We set 
$$P=-t^{d-1}F+t^{d-1}+x(t+xF^d)^d$$
and consider the ideal $\mathcal I\subset \mathcal A_d$
generated by $P_{t^\alpha x^\beta}=
\frac{\partial^{\alpha+\beta}P}{\partial t^\alpha
\partial x^\beta}\in \mathcal A_d$ for $\alpha+\beta\leq d$.

\begin{lem}\label{lemPirred} The polynomial $P=-t^{d-1}F+t^{d-1}+x(t+xF^d)^d$
is irreducible over $\mathbb C[[t,x,F]]$ for every integer $d\geq 1$.
\end{lem}

{\bf Proof} For $d\geq 1$ fixed, consider the Newton-polytope 
$${\mathcal N}(P)=\mathop{Conv}(\lbrace
(a,b,c)\in\mathbb N^3\ \vert (P,t^a x^bF^c)\not=0\rbrace)$$
of $P$
defined as the convex hull of all exponents associated to 
monomials involved in $P$. A straightforward computation shows 
that ${\mathcal N}(P)$ is the $3-$dimensional simplex with vertices
$$(d-1,0,1),(d-1,0,0),(d,1,0),(0,1+d,d^2)\in\mathbb N^3\ .$$
A factorisation $P=P_1P_2$ of $P$ implies
the equality ${\mathcal N}(P)={\mathcal N}(P_1)+{\mathcal N}(P_2)$.
where ${\mathcal N}(P_i)$ is the Newton polytope of the factor $P_i$.
Since ${\mathcal N}(P)$ is a simplex, the polytope ${\mathcal N}(P_i)$ is 
a of the form 
$\lambda_i{\mathcal N}(P)+\tau_i$ with $\lambda_i\in[0,1]$
and $\tau_i\in\mathbb Q^3$. Since the simplex ${\mathcal N}(P)$
has edges without interior integral vertices and since ${\mathcal 
N}(P_i)$ are polytopes with integral vertices, we have
$\{\lambda_1,\lambda_2\}=\{0,1\}$. Suppose $\lambda_1=0$.
The polynomial $P_1$ is thus of the form $\mu t^kx^lW^m$ for some 
$\mu\in\mathbb C^*$ and $(k,l,m)\in\mathbb N^3$
which implies $P_1\in\mathbb C^*$
by inspection of $P$.\hfill$\Box$

We consider now the three elements
$$Q=(D_L-D_R)F,\ \tilde Q=(\tilde D_L-\tilde D_R)F$$ and
$$K=\frac{\partial P}{\partial F}=-t^{d-1}+d^2x^2(t+xF^d)^{d-1}F^{d-1}$$
of $\mathcal A_d$.

\begin{thm}\label{thmequiv} 
Conjecture \ref{conjcoeffwd} holds if and only if we have the 
two inclusions $K^NQ\in\mathcal I$ and
$K^N\tilde Q\in \mathcal I$ for 
$N=\sum_{j=1}^d\lfloor\frac{d}{j}\rfloor$.
\end{thm}

{\bf Proof}
Set $Q_d=Q$. For $h$ such that $1\leq h\leq d$,
define $Q_{h-1}$ by the substitutions
$$F_{t^\alpha x^{h-\alpha}}\longmapsto  
F_{t^\alpha x^{h-\alpha}}-\frac{1}{K}P_{t^\alpha x^{h-\alpha}},\ 0\leq \alpha
\leq h$$
in $K^{d_h}Q_h$ where $d_h=\lfloor d/h\rfloor$ 
is the degree of $Q_h$ with respect to the variables 
$F_{t^\alpha x^{h-\alpha}},\ 0\leq \alpha\leq h$.
One shows by descending induction on $h$ that 
$Q_{h-1}\equiv K^{d_h}Q_h\pmod
{\mathcal I}$ is an element of the algebra $\mathcal A_{h-1}$.
Since $P$ is irreducible by Lemma \ref{lemPirred},
we have $Q_0\in\mathcal I$ if and
only if $Q_0$ is divisible by $P$ which shows that 
$K^NQ\in\mathcal I$ for $N=\sum_{j=1}^d d_j$.

The proof proceeds then as in the case $d=2$. The 
inclusion $K^NQ\in\mathcal I$ (and $K\equiv -t^{d-1}\pmod {x^2}$) implies
$\varphi(Q)=(D_L-D_R)W=0$ where $\varphi:\mathcal A_d
\longrightarrow \mathbf R=\left(\mathbb Q[t,t^{-1}]\right)[[x]]$
is the homomorphism of algebras defined by $\varphi(t)=t,\varphi(x)=x$
and $\varphi(F_{t^\alpha x^\beta})=\frac{\partial^{\alpha+\beta}}{\partial t^\alpha \partial x^\beta}W$ for $W=\sum_{n=0}^\infty w_{d,n}(t)x^n$. 
It contains $\mathcal I$ in
its kernel. Repeating the above arguments with $\tilde Q$ 
ends the proof.\hfill$\Box$

{\bf Proof of Theorem \ref{thm234}} Using Maple 8, we checked the
inclusions $K^NQ,K^N\tilde Q\in \mathcal I$ of Theorem 
\ref{thmequiv} up to $d=6$.
\hfill$\Box$

\begin{rem} The following trick avoids the 
use of rational fractions in the 
substitutions $F_{t^\alpha x^{h-\alpha}}\longmapsto  
F_{t^\alpha x^{h-\alpha}}-\frac{1}{K}P_{t^\alpha x^{h-\alpha}},\ 0\leq \alpha
\leq h$:
Write $Q_h=\sum_{j=0}^{d_h} q_{h,j}$ where $q_{h,j}$ is homogeneous
of degree $j$ with respect to the variables $F_{t^\alpha x^{h-\alpha}}$.
We have then $Q_{h-1}=\sum_{j=0}^{d_h}K^{d_h-j}\tilde q_{h,j}$
where $\tilde q_{h,j}$ is obtained from $q_{h,j}$ by the substitutions
$F_{t^\alpha x^{h-\alpha}}\longmapsto
KF_{t^\alpha x^{h-\alpha}}-P_{t^\alpha x^{h-\alpha}},\ 0\leq \alpha\leq d_h$.

This reduces the computations for proving Theorem \ref{thm234} to 
elementary operations on polynomials, a domain
of excellence for symbolic computer algebra systems.
\end{rem}

{\bf Acknowledgement} I am gratefull to Bernard Parisse for
a Gr\"obner basis computation establishing the very first proof 
of Theorem \ref{thmwindcoeffs} and to Didier Piau for many
suggestions improving the exposition.

%%%%%%%%%%%%%%%%%%%%%%%%%%%%%%%%%%%%%%%%%%%%%%%%%%%%%%%%%%%%%%%%%%%%%%%%%%

\noindent Roland BACHER

\noindent INSTITUT FOURIER

\noindent Laboratoire de Math\'ematiques

\noindent UMR 5582 (UJF-CNRS)

\noindent BP 74

\noindent 38402 St Martin d'H\`eres Cedex (France)
\medskip

\noindent e-mail: Roland.Bacher@ujf-grenoble.fr


\begin{thebibliography}{99}

\bibitem{Ca} D. Callan, Personal communication (e-mail), 17 september 2007.

\bibitem{EIS} N.J.A. Sloane, {\it The On-Line
Encyclopedia of Integer Sequences},
http://www.research.att.com/$\sim$njas/sequences/index.html

\bibitem{HPR} F. Harary, E.M. Palmer, R.C. Read, 
{\it On the cell-growth problem for arbitrary polygons}, Discrete Maths {\bf 11} (1975), 371--389.

\bibitem{Hen} P. Henrici,
{\it Die Lagrange-B\"urmannsche Formel bei formalen
Potenzreihen}, Jahresber. Deutsch. Math.-Verein. {\bf 86} (1984), 115--134).

\bibitem{St2} R.P. Stanley, Enumerative Combinatorics, Volume 2,
Cambridge University Press.

\bibitem{Tak} L. Takacs, {\it Enumeration of rooted trees and forests}, 
Math. Scientist {\bf 18} (1993) No. 1, 1--10.

\end{thebibliography}
\end{document}